\documentclass[11pt]{article}

\usepackage{amssymb}
\usepackage{epsfig}
\usepackage{amsmath}

 \newcommand{\gr}[1]{{{\mbox{\scriptsize $(#1)$}}}}

\newtheorem{theorem}{Theorem}

\newtheorem{corollary}[theorem]{Corollary}
\newtheorem{lemma}[theorem]{Lemma}
\newtheorem{example}[theorem]{Example}
\newcommand{\al}{{\alpha}}
\newcommand{\be}{{\beta}}

\newcommand{\la}{{\lambda}}
\newcommand{\ga}{{\gamma}}
\begin{document}
\centerline{\bf Moduli of Germs of Legendrian Curves }
\bigskip
\centerline{Ant\'onio Ara\'ujo and Orlando Neto}

\begin{abstract}
We construct the generic component of the moduli space
of the germs of Legendrian curves  with generic plane projection 
topologicaly equivalent to a curve $y^n=x^m$.
\end{abstract}

\section{Introduction}
\label{intro}

Zariski \cite{zar} initiated the construction of the
moduli of plane curve singularities. Delorme \cite{del} organized
in a systematic way the ideas of Zariski, obtaining general
results o the case of curves with one characteristic exponent in
the generic case (see also \cite{peraire}). Greuel, Laudal and Pfister (see  the
bibliography of \cite{gre3}) stratified the space versal
deformations of plane curves, constructing moduli spaces on each
stratum. 

In this paper we initiate the study of the moduli of Legendrian
curve singularities. We construct the moduli space of generic irreducible
Legendrian singularities with equisingularity type equal to the topological type 
of the plane curve $y^n=x^m$, $(n,m)=1$. Our method is based on 
the analysis of the action of the
group of infinitesimal contact transformations on the set of Puiseux
expansions of the germs of plane curves.

In section 2 we associate to each pair of positive integers $n,m$ such that $(n,m)=1$
 a semigroup $\Gamma (n,m)$. We show that the semigroup
of a generic element of this equisingularity class equals
$\Gamma (n,m)$. In section 3 we classify the infinitesimal contact
transformations on a contact threefold and study its action on the
Puiseux expansion of a plane curve. In section 4 we discuss some simple examples of
moduli of germs of Legendrian curves. In section 5 we show that the generic components 
of the moduli of germs of Legendrian curves 
with fixed equisingularity
class are the points of a Zariski open subset of a weighted
projective space.

\section{Plane curves versus Legendrian curves}

Let $\Lambda$ be the germ at $o$ of an irreducible space curve. A
local parametrization $\imath: (\mathbb C,0)\to (\Lambda,o)$
defines a morphism $\imath^*$ from the local ring $\mathcal
O_{\Lambda,o}$ into its normalization $\mathbb C\{t\}$.
 The semigroup  of
$\Lambda$ equals the set $\Gamma$ of the orders of the series that
belong to the image of $\imath^*$. There is an integer $k$ such
that $l\in \Gamma$ for all $l\geq k$. The smallest integer $k$
with this property is denoted by $c$ and called the \em conductor
\em of $\Gamma$. 

Let $C$ be the germ at the origin of a singular
irreducible plane curve $C$ parametrized by
\begin{equation}\label{planecurve}
 x=t^n,\qquad y=\sum_{i=m}^\infty a_it^i,
\end{equation}
 with $a_m\not=0$ and $ (n,m)=1$. The pair $(n,m)$ determines
the topological type of $C$.

\bigskip

Let $M$ be a complex manifold of dimension $n$.
The cotangent bundle $\pi_M:T^*M\to M$ of  $M$ is endowed of a canonical $1$-form $\theta$.
The differential form $(d\theta)^{\wedge n}$ never vanishes on $M$. Hence
$d\theta$ is a symplectic form on $T^*M$. Given a system of local coordinates $(x_1,\ldots,x_n)$ 
on an open set $U$ of $X$, there are holomorphic functions $\xi_1,\ldots,\xi_n$ 
on $\pi_M^{-1}(U)$ such that $\theta\mid_{\pi_M^{-1}(U)}=\xi_1dx_1+\cdots+\xi_ndx_n$.

Let $X$ be a complex threefold.
Let $\Omega^k_X$ denote the sheaf of differential forms of degree $k$ on $X$.
A local section of $\Omega^1_X$ is called a \em contact form \em if $\omega \wedge d\omega$ never vanishes.
Let $\mathcal L$ be a subsheaf of the sheaf $\Omega_X^1$.
The sheaf $\mathcal L$ is called a \em contact structure \em on $X$ if $\mathcal L$ is locally generated by a contact form.
A pair $(X,\mathcal L)$, where $\mathcal L$ is a contact structure on $X$, is called a \em contact threefold. \em
Let $(X_i,\mathcal L_i)$, $i=1,2$, be two contact threefolds. A holomorphic map $\varphi:X_1\to X_2$ 
is called a \em contact transformation \em if $\varphi^*{\mathcal L}_2=\mathcal L_1$.

Let $\mathbb{P}^{*}\mathbb{C}^2=\mathbb{C}^2\times
\mathbb{P}^1=\{(x,y,(\xi :\eta )):x,y,\xi ,\eta \in \mathbb{C}$,
$(\xi ,\eta )\neq (0,0)\}$ be the projective cotangent bundle of
$\mathbb{C}^2$. Let $\pi :\mathbb{P}^{*}\mathbb{C}^2\rightarrow
\mathbb{C}^2$ be the canonical projection. Let $U$ and $V$ be the
open sets of $ \mathbb{P}^{*}\mathbb{C}^2$ defined respectively by
$\eta \neq 0$ and $\xi \neq 0$. Set $p=-\xi /\eta $, $q=-\eta /\xi
$. The sheaf $\mathcal{L}$ defined by $\mathcal{L\mid
}_U=\mathcal{O}_U(dy-pdx)$ and $\mathcal{L\mid
}_V=\mathcal{O}_V(dx-qdy)$ is a contact structure on $
\mathbb{P}^{*}\mathbb{C}^2.$ By the Darboux theorem every contact
threefold is locally isomorphic to $(U,\mathcal{O}_U(dy-pdx))$. We
call \em infinitesimal contact transformation \em to a germ of a
contact transformation $\Phi:(U,0) \mapsto (U,0)$.

A curve $\Lambda $ on a contact manifold $(X,\mathcal{L})$ is
called \em Legendrian \em if the restriction of $\omega $ to the
regular part of $\Lambda $ vanishes for each section $\omega $ of
$\mathcal{L}$. Let $C=\{f=0\}$ be a plane curve. Let $\Lambda $ be
the closure on $\mathbb{P} ^{*}\mathbb{C}^2$ of the graph of the
Gauss map $G:\{a\in C:df(a)\neq 0\}\rightarrow \mathbb{P}^1$
defined by $G(a)=\langle df(a)\rangle $. The set $\Lambda $ is a
Legendrian curve. We call $\Lambda $ the \em conormal \em of the
curve $C.$ If $C$ is irreducible and parametrized by
(\ref{planecurve}) then $\Lambda $ is parametrized by
\begin{equation}\label{legc}
x=t^n,\qquad y=\sum_{i=m}^\infty a_it^i, \qquad
p=\frac{dy}{dx}=\sum_{i=m}^\infty \frac{i}{n} a_i t^{i-n}.
\end{equation}
Given a Legendrian curve $\Lambda $ of
$\mathbb{P}^{*}\mathbb{C}^2$ such that $ \Lambda $ does not
contain any fibre of $\pi $, $\pi (\Lambda )$ is a plane curve.
Moreover, $\Lambda $ equals the conormal of $\pi (\Lambda )$.

Let $(X,\mathcal{L})$ be a contact threefold. A holomorphic map
$\varphi :(X,o)\rightarrow (\mathbb{C}^2,0)$ is called a \em
Legendrian map \em if $ D\varphi (o)$ is surjective and the fibers
of $\varphi $ are smooth Legendrian curves.\thinspace The map
$\varphi \,$ is Legendrian if and only if there is a contact
transformation $\psi :(X,o)\rightarrow (\mathbb{P}^{*}
\mathbb{C}^2,(0,0,(0:1))$ such that $\varphi =\pi \psi .$

 Let $(\Lambda ,o)$ be a
Legendrian curve of $ X$. Let $C_o(\Lambda )$ be the tangent cone
of $\Lambda $ at $o$. We say that a Legendrian map $\varphi
:(X,o)\rightarrow (\mathbb{C}^2,0)$ is \em generic \em relatively
to $(\Lambda ,o)$ if it verifies the transversality condition $
T_o\varphi ^{-1}(0)\cap C_o(\Lambda )=\{0\}.$ We say that a
Legendrian curve $(\Lambda ,o)$ of $\mathbb{P}^{*}\mathbb{C}^2$ is
in strong generic position if $\pi
:(\mathbb{P}^{*}\mathbb{C}^2,o)\rightarrow (\mathbb{C} ^2,\pi
(o))$ is generic relatively to $(\Lambda ,o)$. The Legendrian
curve $ \Lambda $ parametrized by (\ref{legc}) is in \em strong
generic position \em if and only if $ m\geq 2n+1.$ Given a
Legendrian curve $(\Lambda ,o)$ of a contact threefold $X$ there
is a contact transformation $\psi :(X,o)\rightarrow
(\mathbb{P}^{*}\mathbb{C} ^2,(0,0,(0:1))$ such that $(\psi
(\Lambda ),o)$ is in strong generic position (cf \cite{kk},
section 1).

We say that two germs of Legendrian curves are
equisingular if their images by generic Legendrian maps have 
the same topological type.

\section{Infinitesimal Contact Transformations}
Let $m$ be the maximal ideal of the ring $\mathbb{C}\{x,y,p\}$.
 Let $\mathcal{G}$ denote the group of infinitesimal
contact transformations $\Phi$ such that the derivative of $\Phi$
leaves invariant the tangent space at the origin of the curve
$\{y=p=0\}$. Let $\mathcal{J}$ be the group of infinitesimal
contact transformations
\begin{equation}\label{ctransf} (x,y,p)\mapsto (x+\alpha ,y+\beta
,p+\gamma )
\end{equation}
 such that $\al,\be,\ga,\partial\al /\partial x,\partial\be
 /\partial y, \partial\ga /\partial p\in m$. Set
 $\mathcal{H}=\{\Psi_{\lambda , \mu}:\la,\mu\in \mathbb{C}^*\}$, where
\begin{equation}\label{homo}
\Psi_{\la,\mu}(x,y,p) = \left(\lambda x,\mu y, \frac{\mu }{\lambda
}p\right).
\end{equation}
 Let $\mathcal{P}$ denote the group of
\em paraboloidal contact transformations \em (see \cite{k2})
\begin{equation}\label{para}
(x,y,p) \mapsto (ax+bp,y- \frac{1}{2} acx^2- \frac{1}{2}
bdp^2-bcxp,cx+dp),~\left|
\begin{array}{ll}
a & b \\
c & d \\
\end{array}
\right|=1.
\end{equation}
The contact transformation (\ref{para}) belongs to $\mathcal{G}$
if and only if $c=0$. The paraboloidal contact transformation
\begin{equation}\label{paraleg}
(x,y,p) \mapsto (-p,y-xp,x)
\end{equation}
Is called the \em Legendre transformation.\em
\begin{theorem}\label{subgroup}
The group $\mathcal{J}$ is an invariant subgroup of $\mathcal{G}$.
Moreover, the quotient $\mathcal{G}/\mathcal{J}$ is isomorphic to
$\mathcal{H}$.
\end{theorem}
Proof.  If $H\in \mathcal{H}$  and $\Phi\in \mathcal{J}$,  $H \Phi H^{-1} \in
\mathcal{J}$. Hence it is enough to show that each element of
$\mathcal{G}$ is a composition of  elements of $\mathcal{H}$ and
$\mathcal{J}$. Let $\Phi\in\mathcal{G}$ be the infinitesimal
contact transformation $(x,y,p) \mapsto (x',y',p')$. There is
$\varphi \in \mathbb{C}\{x,y,p\}$ such that $\varphi (0) \neq 0$
and
\begin{equation}\label{eqq}
  dy'-p'dx'=\varphi(dy-pdx).
\end{equation}
Composing $\Phi$  with $H\in\mathcal H$ we can assume that
$\varphi (0)=1$. Let $\hat{\Phi}$ be the germ of the symplectic
transformation $(x,y,p;\eta)\mapsto (x',y',-\eta
p';\varphi^{-1}\eta )$. Notice that
$\hat{\Phi}(0,0;0,1)=(0,0;0,1)$. Since $D \hat{ \Phi}(0,0;0,1)$
leaves invariant the linear subspace $\mu$ generated by
$(0,0;0,1)$, $D \hat{ \Phi}(0,0;0,1)$ induces a linear symplectic
transformation on the linear symplectic space $\mu ^\bot / \mu$.
There is a paraboloidal contact transformation $P$ such that
$D\hat{P}(0,0;0,1)$ equals $D \hat{ \Phi}(0,0;0,1)$
 on $\mu ^\bot / \mu$. Since $D(\hat P^{-1}\hat\Phi)(0,0;0,1)$
 induces the identity map on $\mu ^\bot / \mu$,
$P^{-1}\Phi$ is an infinitesimal contact transformation of the
type $ (x,y,p) \mapsto (x+\alpha, y', p+\gamma)$, where
\begin{equation}\label{eqqbis}
  {\partial \alpha \over \partial x},{\partial \alpha \over
\partial p},{\partial \gamma \over \partial x},{\partial \gamma
\over \partial p} \in m.
\end{equation}
Set $\beta=y'-y$. It follows from (\ref{eqq}) and (\ref{eqqbis})
that $(\partial\be /\partial y)(0)=0$. Hence $P^{-1}\Phi \in
\mathcal{J}$. Since $\Phi$ and $P^{-1}\Phi \in \mathcal{G}$,
$P\in\mathcal{G}$. Therefore $p$ is the composition of an element
of $\mathcal{H}$ and an element of $\mathcal{J}$.
%
q.e.d.

\begin{theorem}\label{existence}
Let $\alpha \in \mathbb{C}\{x,y,p\},\beta _{0}\in
\mathbb{C}\{x,y\}$ be power series such that
\begin{equation}\label{alfam}
\alpha ,\beta _{0},\frac{\partial \beta _{0}}{\partial y}\in m.
\end{equation}
\em There are $\beta ,\gamma  \in \mathbb{C}\{x,y,p\}$ such that $
\beta -\beta _{0}\in (p)$, $\gamma\in m $ and $\alpha,\beta,\gamma$ define an infinitesimal
 contact transformation $\Phi_{\alpha,\beta_0}$ of type (\ref{ctransf}). The power series $\beta$
and $\gamma$ are uniquely determined by these conditions.
Moreover, $(\ref{ctransf})$ belongs to $\mathcal{J}$ if and only if
\begin{equation}\label{alfama}
\frac{\partial \alpha}{\partial x}, \frac{\partial \beta
_{0}}{\partial x},\frac{\partial^2 \beta _{0}}{\partial x\partial
p}\in m.
\end{equation}
The function $\beta$ is the solution of the Cauchy problem
\begin{equation}\label{CK}
\left(1+\frac{\partial \alpha }{\partial x}
+p\frac{\partial \alpha }{\partial y}\right)
\frac{\partial \beta }{\partial p}-p\frac{\partial \alpha
}{\partial p}\frac{\partial \beta }{\partial y}
-\frac{\partial \alpha }{\partial p}\frac{
\partial \beta }{\partial x} =p\frac{\partial \alpha }{\partial p}.
\end{equation}
with initial condition $\beta-\beta_0\in (p)$.

\end{theorem}
Proof. The map (\ref{ctransf}) is a contact
transformation if and only if there is $\varphi \in
\mathbb{C}\{x,y,p\}$ such that $\varphi (0) \neq 0$ and
\begin{equation}\label{cont0}
d(y+\beta )-(p+\gamma )d(x+\alpha )=\varphi (dy-pdx).
\end{equation}
The equation (\ref{cont0}) is equivalent to the system
\begin{eqnarray}\label{sist1}
 \frac{\partial \beta }{\partial p} & = &
(p+\gamma )\frac{\partial\alpha }{\partial p} \\
\label{sist2}
 \varphi & = & 1+\frac{\partial \beta }{\partial
y}-(p+\gamma )\frac{\partial \alpha }{\partial y} \\
\label{sist3}
 -p\varphi & = & \frac{\partial \beta }{\partial
x}-(p+\gamma )(1+\frac{\partial \alpha }{\partial x}).
 \end{eqnarray}
By (\ref{sist2}) and (\ref{sist3}),
\begin{equation}\label{CKga}
\frac{\partial \be }{\partial x}
-(p+\ga )\left(1+\frac{\partial \al}{\partial x}
+p\frac{\partial\al }{\partial y}   \right)
+p\left(1+\frac{\partial \be }{\partial y}
\right)=  0,
\end{equation}
By (\ref{sist1}) and (\ref{CKga}), (\ref{CK}) holds.

By the Cauchy-Kowalevsky theorem there is one and only one
solution $\beta $ of (\ref{CK}) such that $\beta -\beta _{0}\in
(p)$. It follows from (\ref{CKga}) that
\begin{equation}\label{gama}
\gamma =\left(1+\frac{\partial \alpha }{\partial
x}+p\frac{\partial \alpha }{
\partial y}\right)^{-1}
\left(\frac{\partial \beta }{\partial x}+
p\left(\frac{\partial \beta }{\partial y}-\frac{\partial \alpha
}{\partial x}-p\frac{\partial \alpha }{
\partial y}\right)\right).
\end{equation}
Since $\partial \beta _{0}/\partial y\in m$, $\partial \beta
/\partial y\in m$. By (\ref{sist2}), $\varphi (0)\neq 0.$

\noindent (ii)  Since $\partial\be_0 /\partial x\in m$,
$\partial\be /\partial x\in m$. By (\ref{gama}), $\ga\in m$. By
(\ref{gama}),
\[
\frac{\partial\ga}{\partial p}\in
\left( \frac{\partial^2\be}{\partial x\partial
p}+\frac{\partial\be}{\partial y}-\frac{\partial\al}{\partial x},p
\right).
\]
By (\ref{alfam}) and (\ref{alfama}), $\partial\ga /\partial p\in
m$.
q.e.d.
%
%
%

\begin{corollary}\label{group}
The elements of $\mathcal J$ are the infinitesimal contact transformations $\Phi_{\alpha,\beta_0}$ such that 
 $\alpha,\beta _{0}$ verify $($\ref{alfam}$)$ and $($\ref{alfama}$)$.
\end{corollary}

\begin{lemma}\label{forget}
Given $\lambda\in\mathbb C$ and $w\in \Gamma (m,n)$ such that
$w\ge m+n$, there are $\alpha,\beta_0$ verifying the conditions of
theorem \ref{existence} such that $\imath^*(\beta-p\alpha)=\lambda
t^w+\cdots$.
\end{lemma}
Proof. 
By (\ref{estimation}) there is $b\in\mathbb C\{x,y,p\}$ such that
$\imath^*b=\lambda t^w+\cdots$, $b=\sum_{k\ge 0}b_kp^k$ and
$v(b_k)\ge v(b)-v(x)-kv(p)+1$. Set $\alpha=-\partial b/\partial p$,
$\beta_0=b_0$. Set $\al=\sum_{k\ge 0}\al_kp^k$, $\be=\sum_{k\ge
0}\be_kp^k$, where $\al_k,\be_k\in\mathbb{C}\{x,y\}$.
By (\ref{CK}),
\begin{eqnarray*}
k\be_k   &  +   &   \sum_{j=1}^{k-1} j\be_j
\left( \frac{\partial \al_{k-j}}{\partial x} +
\frac{\partial \al _{k-j-1} }{\partial y} \right)= \cr
&  =  &  (k-1)\al_{k-1}+k\al_k\frac{\partial \be_0}{\partial x}+
\sum_{j=1}^{k-1}j\al_j\left(\frac{\partial \be_{k-j}}{\partial x}
+\frac{\partial \be_{k-j-1}}{\partial y}\right),
\end{eqnarray*}
for $k\ge 1$. Since $\al_l=-(l+1)b_{l+1}$ for $l\ge 1$, $ v(\al_jp^k)\ge w+1$,
if $j\le k-2$. Moreover, $v(\al_{k-1}p^k)\ge w+1-n$ and $v(\al_{k}p^k)\ge w+1-m$.
Therefore
\[
k\be_k p^k   +     \sum_{j=1}^{k-1} j\be_j
\left( \frac{\partial \al_{k-j}}{\partial x}  +
\frac{\partial \al _{k-j-1} }{\partial y} \right)p^k \equiv
(k-1)\al_{k-1}p^k  + (k-1)\al_{k-1}\frac{\partial \beta_1}{\partial x} ,
\]
mod $\left( t^{w+1} \right)$ for $k\ge 1$. We show by induction in $k$ that 
\[ k\be_kp^k\equiv
(k-1)\al_{k-1}p^k\  \mod ( t^{w+1} ) ,\hbox{\rm  for } k\ge 1.
\]
 Hence
 $\beta-p\al\equiv b$ mod $( t^{w+1} )$.
%
%
q.e.d.
There is an action of $\mathcal J$ into the set of germs  of plane curves
$C$ such that the tangent cone to the conormal of $C$ 
equals $\{y=p=0\}$.  Given $\Phi\in \mathcal J$ we associate to $C$ 
the image by $\pi\Phi$ of the conormal of $C$. Given integers $n,m$
such that $ (m,n)=1$ and $m\ge 2n+1$, $\mathcal J$ acts on the
series of type (\ref{planecurve}). Given an infinitesimal contact transformation (\ref{ctransf}) there is $s \in
\mathbb{C}\{t\}$ such that $s^n=t^n+\alpha$ and for each $i \geq
1$
\[
s^i=t^i\left(1+\frac{i}{n}\frac{\alpha (t)
}{t^n}+\frac{i}{n}\left(\frac{i}{n}-1\right)\left(\frac{\alpha (t)
}{t^n}\right)^2+\cdots\right).
\]

\begin{lemma}\label{action} If $v(\beta_0)\ge v(\alpha )+v(p)$,
 the contact transformation \em (\ref{ctransf}) \em takes
  \em (\ref{planecurve}) \em into the plane curve parametrized
by $ x =s^n$, $y=y(s) + \beta(s) -p(s)\alpha(s) + \varepsilon $,
where $v(\varepsilon) \geq 2v(\alpha)
 + m - 2n. $
\end{lemma}
Proof. Since
$t^i=s^i-(i/n)t^{i-n}\alpha (t) +(i(i-n)/n^2)\alpha (t)^2
t^{i-2n}+\cdots$,
\[
y(t)=\sum_{i \geq m}a_i s^i - \alpha (t) \sum_{i \geq m}
\frac{i}{n} a_i t^{i-m} + \varepsilon'=y(s)-\alpha
(t)p(t)+\varepsilon',
\]
\[
p(t)\alpha (t)=p(s)\alpha (t)-\alpha (t)^2\sum_{i \geq m}
(\frac{i}{n})^2 a_i t^{i-2m} + \varepsilon''=p(s)\alpha
(s)+\varepsilon''',
\]
where $v(\varepsilon' ),v(\varepsilon'' ),v(\varepsilon''' )\ge
2v(\al)+m-2n $.
q.e.d.


\section{Examples}

\begin{example}
If $m$ odd all plane curves topologicaly equivalent to $y^2=x^m$ are analyticaly equivalent to $y^2=x^m$ $($cf. \cite{zar}$)$. 
Hence all Legendrian curves with generical plane projection $y^2=x^m$ are contact equivalent to the conormal of $y^2=x^m$.
\end{example}

\begin{example}
Let $m,s,\epsilon$ be positive integers. Assume that $m=3s+\epsilon$, $1\le \epsilon\le 2$. 
Let $C_{3,m,\nu}$ be the plane curve parametrized by 
\[
x=t^3, ~~~ y=t^m+t^{m+3\nu+\epsilon-3}.
\]
By \cite{zar} a plane curve topologically equivalent to $y^3=x^m$ is analyticaly equivalent to $y^3=x^m$ 
or to one of the curves $C_{3,m,\nu}$, $1\le \nu\le s-1$. The infinitesimal contact transformation
\[
(x,y,p)\mapsto (x-2p,y+p^2,p)
\]
takes the plane curve $C_{3,m,s-1}$ into the plane curve $C'$ parametrized by
\[
3x=3t^3-mt^{m-3}-\cdots, ~~~ y=t^m.
\]
By Lemma \ref{action}, the curve $C'$ admits a parametrization of the type $x=s^3$, $y=s^m+\delta$, where $v(\delta)\ge m+3s+\epsilon-6$.
By \cite{zar}, the curve $C'$ is analyticaly equivalent to the plane curve $y^3=x^m$.  

The semigroup of the conormal of the plane curve $y^3=x^m$ equals $\Gamma _{3,m,0}=\langle 3,m-3\rangle$. 
The semigroup of the conormal of the curve $C_{3,m,\nu}$ equals $\Gamma_{3,m,\nu}=\langle 3,m-3, m+3\nu+\epsilon\rangle$, $1\le \nu\le s-1$. 
The map from $\{0,1,\ldots,s-2\}$ into $\mathcal P(\mathbf N)$ that takes $\nu$ into $\Gamma_{3,m,\nu}$ is injective. 
Hence there are $s-1$ analytic equivalence classes of plane curves topologicaly equivalent to $y^3=x^m$ and 
$s-2$ equivalence contact classes of Legendrian curves with generical plane projection $y^3=x^m$. In this case the semigroup of a
 curve is an analytic invariant that classifies the contact equivalence classes of Legendrian curves. 
We will see that in the general case there are no discrete invariants that can classify the contact equivalence classes of Legendrian curves.
\end{example}
Given a plane curve 
\begin{equation}\label{genn}
x=t^3, ~~~ y=t^m+\sum_{i\ge m+\epsilon}a_itî,
\end{equation}
the semigroup of the conormal of (\ref{genn})  equals $\Gamma_{3,m,1}$
if and only if $a_{m+\epsilon}\not=0$. It is therefore natural to call $\Gamma (3,m):=\Gamma_{3,m,1}$ 
the generic semigroup of the family of Legendrian curves with generic plane projection $y^3=x^m$.

\section{The generic semigroup of an equisingularity class of irreducible
Legendrian curves}

We will  associate to a pair $(n,m)$ such that $m\ge 2n+1$ and
$(m,n)=1$ a semigroup $\Gamma (n,m)$. Let $\langle
k_1,\ldots,k_r\rangle$ be the submonoid of $\mathbb (N,+)$
generated by $k_1,\ldots,k_r$. Let $c$ be the conductor of the
semigroup of the plane curve (\ref{planecurve}). Set $\Gamma
_c=\langle n\rangle\cup \{c,c+1,...\}.$ We say that the trajectory
of $k\ge c$ equals $\{k,k+1,...\}$.
 Let
us assume that we have defined $\Gamma _j$ and the trajectory of
$j$ for some $j\in \langle n,m-n\rangle\setminus\Gamma_c$, $j\ge
m$. Let $i$ be the biggest element of $\langle
n,m-n\rangle\setminus \Gamma _j$. Let $\sharp_{i}$ be the minimum
of the cardinality of the set of monomials of $\mathbb{C[} x,y,p]$
of valuation $i$ and the cardinality of $\{ i,i+1,\ldots
\}\backslash \Gamma_j$. Let $\omega_i$ be the $\sharp_{i}$-th
element of $\{ i,i+1,\ldots \}\backslash \Gamma_j.$ We call \em
trajectory \em of $i$ to the set $\tau_i = \{ i,i+1,\ldots,
\omega_i\}\setminus \langle n\rangle$. Set $\Gamma_{i}= \tau_{i}
\bigcup \Gamma_j$. Set $\Gamma (n,m)=\Gamma_{m-n}$. 
The main
purpose of this section is to prove theorem \ref{deter}. Let us
show that
\begin{equation}\label{estimation}
  \omega_i\le i+n-2.
\end{equation}
If $\omega_i\ge i+n-1$ , $\Gamma_i\supset \{i,\ldots,i+n-1\}$. 
Hence $\Gamma_{i}\supset \{i,i+1,\ldots\}$ and $i\ge c$. Therefore
(\ref{estimation}) holds.

Let $X=t^n$, $ Y=\sum_{i\ge 0}a_{m+i}t^{m+i}$, $P=\sum_{i\ge
0}(\mu+i)a_{m+i}t^{m-n+i}$  be power series with coefficients in
the ring $\mathbb{Z}[a_m,\ldots,a_{c-1},\mu]$. Given $J=(i,j,l)\in
\mathbb N^3$, set $v(J)=v(x^iy^jp^l)$. Let $\mathcal N=\{J\in
\mathbb N^3 : j+l\ge 1$ and $v(J)\le c-1\}$.
Let $\Upsilon=(\Upsilon_{J,k})$, $J\in \mathcal N$, $m\le k\le
c-1$ be the matrix such that
\begin{equation}\label{defupsilon}
X^iY^{j}P^l\equiv\sum_{k= m}^{c-1}\Upsilon_{J,k}t^{k} \qquad (mod
~(t^c)).
\end{equation}
Since $\partial Y /
\partial \mu =0$ and $X\partial P /
\partial \mu =Y$,
\begin{equation}\label{derivadaYP}
\frac{\partial X^iY^j P^l}{\partial
\mu}=lX^{i-1}Y^{j+1}P^{l-1}\qquad \hbox{\rm and} \qquad
\frac{\partial \Upsilon_{J,k}}{\partial \mu} = l
\Upsilon_{\partial J,k},
\end{equation}
where $\partial (i+1,j,l+1)=(i,j+1,l)$. Moreover,
\begin{equation}\label{YNPJmatrix}
\Upsilon_{J,k}=\sum_{\alpha \in A (k) }\sum_{\gamma \in G
(\alpha,l)}\frac{j!~l!}{(\alpha-\gamma)!\gamma!}~a^{\alpha} \mu
^{\gamma},
\end{equation}
where $A(k)  =  \{\alpha=(\alpha_m,...,\alpha_{c-1}) :
|\alpha|=j+l$ and $\sum_{s=m}^{c-1} s \alpha_s=k-(i-l)n\}$, $
G(\alpha,l) = \{\gamma :|\gamma|=l $ and $0\le\gamma\leq\alpha \}$
and $ \mu ^{\gamma} = \prod^{c-1}_{s=m}(\mu-m+s)^{\gamma_s}$. Let
us prove (\ref{YNPJmatrix}). We can assume that $i=l$. Since
$G(\alpha,N)=\{\alpha\}$  and $X^NP^N= \sum_{k\ge 0
}t^{k}\sum_{\alpha\in A(k)}
 (N!/\alpha!)\mu^\alpha a^\alpha$ ,
(\ref{YNPJmatrix}) holds for $J=(N,0,N)$. Let us show by induction
in $j$ that (\ref{YNPJmatrix}) holds when $j+l=N$.
 Set $e_s=(\delta_{s,r})$, $0\le s,r\le
N$. Given $\gamma\in G(\alpha,l-1)$, set
$\gamma_{(s)}=\gamma+e_s$. Set $\Delta^\gamma_s=1$ if
$\gamma_{(s)}\le \alpha$. Otherwise, set $\Delta^\gamma_s=0$.
Since
\begin{eqnarray*}
\frac{1}{l}\sum_{\gamma \in G
(\alpha,l)}\frac{j!l!}{(\alpha-\gamma)!\gamma!}\frac{\partial \mu
^{\gamma}}{\partial \mu}
 & = & \sum_{\gamma \in G
(\alpha,l-1)} \sum_{s=m}^{c-1}
\frac{j!(l-1)!}{(\alpha-\gamma_{(s)})!\gamma_{(s)}!}(\gamma_s+1)\Delta_s^\gamma\mu^\gamma
  \\
& = & \sum_{\gamma \in G (\alpha,l-1)}
\frac{j!(l-1)!}{(\alpha-\gamma)!\gamma!}     \mu^\gamma
\sum_{s=m}^{c-1}(\alpha_s-\gamma_s)
  \\
& = & \sum_{\gamma \in G
(\alpha,l-1)}\frac{(j+1)!(l-1)!}{(\alpha-\gamma)!\gamma!}
\mu^\gamma,
\end{eqnarray*}
the induction step follows from (\ref{derivadaYP}).
We will consider in the polynomial ring $\mathbb
C[a_m,\ldots,a_{c-1} ]$ the order $a^\alpha< a^\beta$ 
if there is an integer $q$ such that 
$\alpha_q<\beta_q$ and $\alpha_i=\beta_i$ for $i\ge q+1$. Set
$\omega (P)=\sup\{i:a_i $ occurs  in $P\}$.
%
%
%

\begin{lemma}\label{DDD} Let $M,N,q\in \mathbb Z$ such that
$0\le M\le N$ and $q+N\ge 0$. If $\lambda=(\lambda_{l,k})$, where
$M\le l\le N$, $k\ge 0$, $\lambda_{l,k}=\Upsilon_{J,k}$ and
$J=(q+l,N-l,l)$,  the minors of $\lambda$ with $N-M+1$ columns
different from zero do not vanish at $\mu=m$.
\end{lemma}
Proof. One can assume that $q=0$. When we multiply the left-hand side of
(\ref{defupsilon}) by $P$ the coefficients of $\Upsilon$ are
shifted and multiplied by an invertible matrix. Hence one can
assume that $M=0$. Set $Z=(Z_{j,k})$, where $Z_{j,k}=
{j\choose
k}\mu^{j-k}$, $0\le j,k\le N$. Notice that $Z$ is lower diagonal,
$\det(Z)=1$ and
\begin{equation}\label{ZZZ}
\frac{\partial Z_{j,k} }{ \partial \mu } ~ = ~ jZ_{j-1,k} ~ = ~
(k+1) Z_{j,k+1}.
\end{equation}
Let us show that
\begin{equation}\label{DET}
Z^{-1}\lambda=\lambda|_{\mu=0}.
\end{equation}
Since $\lambda_{N,k}$ is a polynomial of degree $N$ in the
variable $\mu$ with coefficients in the ring $\mathbb Z
[a_m,\ldots,a_{c-1}]$, there are polynomials $\mathcal
Z_{i,k}\in\mathbb Q [a_m,\ldots,a_{c-1}]$ such that
$\lambda_{N,k}=\sum_{i=0}^N{N\choose i}\mathcal Z_{i,k}\mu^{N-i}$.
Set $\mathcal Z=\left( \mathcal Z_{i,k}\right)$, $0\le i\le N$,
$0\le k\le c-1$. Since $Z|_{\mu=0}=Id$, it is enough to show that
$Z\mathcal Z=\lambda$. By construction,
\begin{equation}\label{evaluation2}
\lambda_{j,k}=\sum_{i=0}^NZ_{j,i}\mathcal Z_{i,k}
\end{equation}
when $j=N$.  By (\ref{derivadaYP}) and (\ref{ZZZ}) statement
 (\ref{evaluation2}) holds for all $j$.
Remark that
\begin{equation}\label{e1}
\lambda_{l,v(J)+k}|_{\mu=0}=0\qquad \hbox{\rm if and only
if}\qquad k<l.
\end{equation}
Let $\theta_{l,k}$ be the leading monomial of $\lambda_{l,k}$.
When $k\ge l$,
\begin{equation}\label{e2}
\theta_{l,v(J)+k}=a_m^{N-1}a_{m+k}\qquad \hbox{\rm if} \qquad l=0,
\end{equation}
\begin{equation}\label{e3}
\theta_{l,v(J)+k}= a_m^{N-l}a_{m+1}^{l-1}a_{m+k-l+1}\qquad
\hbox{\rm if} \qquad l\ge 1.
\end{equation}
Let us prove (\ref{e3}).
 Set $\alpha_0=j$, $\alpha_1=l-1$, $\alpha_{k-l+1}=1$ and
$\alpha_s=0$ otherwise. By (\ref{YNPJmatrix}), $\alpha\in A(k)$
and there is one and only one $\gamma\in G(\alpha,j)$ such that
$\gamma_0=0$, the tuple $\overline{\alpha}$ given by
$\overline{\alpha}_0=0$ and $\overline{\alpha}_i=\alpha_i$ if
$i\not =0$. Since
\[
\sum_{\gamma \in G(\alpha,l)} \frac{j!l!\mu
^{\gamma}}{(\alpha-\gamma)!\gamma!}
\equiv
\frac{j!l!\mu^{\overline{\alpha}}}{(\alpha-\overline{\alpha})!\overline{\alpha}!}
\equiv l \prod_{s=0}^{c-m-1} s^{{\overline{\alpha}}_s}= (k-l+1)l
\mod \mu,
\]
the coefficient of $a_m^{N-l}a_{m+1}^{l-1}a_{k-l+1}$ does not
vanish. By (\ref{YNPJmatrix}), $\alpha_{k-l+r} \neq 0$ for some
$r>1$ implies that $\gamma_0>0$ for all $\gamma \in G(\alpha,l)$.
Hence (\ref{e3}) holds.

Let $\lambda'$  be the square submatrix of $\lambda$ with columns
$g\gr{i}+Nm$,
$0  \le  g
\gr{0} <   \cdots   <  g \gr{N} $.
 By (\ref{DET}), $\det (\lambda'|_{\mu=0})=\det
(Z^{-1}\lambda')=\det (Z)^{-1}\det \lambda'=\det \lambda'$. Hence
$\det\lambda'$ does not depend on $\mu$ and $\det
(\lambda'|_{\mu=m})=\det (\lambda'|_{\mu=0})$. Set $
\det(\lambda')=\sum_{\pi}sgn(\pi)\lambda_{\pi}$, where
$\lambda_{\pi}=\prod_{i=0}^N \lambda'_{i,\pi(i)}$.  If
$\lambda_\pi\not=0$, let $\theta_\pi$ be the leading monomial  of
$\lambda_\pi$.

Let $\varepsilon$ be the following permutation of
$\{0,\ldots,N\}$. Assume that $\varepsilon $ is defined for $0\le
i\le l-1$. Let $p_l$ and $q_l$ be respectively the maximum and the
minimum  of $\{0,\ldots,N\}\setminus \varepsilon
\gr{\{0,\ldots,l-1\}}$. If $\lambda_{l+1,q_l}=0$, set $\varepsilon
\gr{l}=q_l$. Otherwise, set $\varepsilon \gr{l}=p_l$. Let us show
that (\ref{e1}) implies that $\lambda_\varepsilon\not=0$. It is
enough to show that $\lambda_{i,q_i}\not=0$ for all $i$. Since
$g\gr{0}\ge 0$, $\lambda_{0,q_0}\not=0$. Assume that $l\ge 1$ and
$\lambda_{i,q_i}\not=0$ for  $0\le i\le l-1$. Hence
$g\gr{q_{l-1}}\ge l-1$.
If $\lambda_{l,q_{l-1}}\not=0$ then $\lambda_{l,q_{l}}\not=0$. If
$\lambda_{l,q_{l-1}}=0$ then $\varepsilon\gr{l-1}=q_{l-1}$.
Therefore $g\gr{q_{l}}=g\gr{q_{l-1}+1}\ge g\gr{q_{l-1}}+1\ge l$
and $\lambda _{l,q_l}\not= 0$.

Let us show that $\theta_\varepsilon$ is the leading monomial of
$\det (\lambda'|_{\mu=0})$. Let $\pi$ be a permutation of
$\{1,\ldots,N\}$.  Assume that $\pi(i)=\varepsilon (i)$ if $0\le
i\le l-1$ and $\pi (l)\not=\varepsilon (l)$. If
$\lambda_{l,q_{l-1}}=0$ then $\pi (l)\not=q_l$ and
$\lambda_\pi=0$. If $\lambda_{l,q_{l-1}}\not=0$  then $\pi
(l)\not=p_l$ and $\omega (\prod_{i=l}^N\lambda_{i,\pi (i)})<
\omega (\prod_{i=l}^N\lambda_{i,\varepsilon (i)})$. Therefore
$\lambda_\pi<\lambda_\varepsilon$.
q.e.d.
The semigroup of the legendrian curve (\ref{legc}) only depends on
$(a_m,\ldots,a_{c-1})$. We will denote it by
$\Gamma_{(a_m,\ldots,a_{c-1})}$.

\begin{theorem}\label{deter}
There is a dense Zariski open subset $U$ of $\mathbb C^{c-m}$ such
that if $(a_m,\ldots,a_{c-1})\in U$,
$\Gamma_{(a_m,\ldots,a_{c-1})}=\Gamma (n,m)$.
\end{theorem}
Proof. 
Since $U$ is defined by the non vanishing of several determinants,
it is enough to show that $U\not=\emptyset$. Let $j\in\langle
n,m-n\rangle$, $j\ge m$. Set $q=\sharp(\tau_j)$. Assume that we
associate to $j$ a family of triples $I_1,\ldots,I_q\in\mathcal N$
such that $v(I_s)\ge j$, $1\le s\le q$, and if $E$ is the linear
subspace of $\mathbb C[a_m,\ldots,a_{c-1}]\{t\}$ spanned by
$\Upsilon_{I_s,k}|_{\mu=m}$, $1\le s\le q$,
$v(E)=\tau_j\cup\{\infty\}$. Let $i$ be the biggest element of
$\langle n,m-n\rangle\setminus\Gamma_j$. Assume that
$\tau_i\cap\tau_j\not=\emptyset$. Hence $\tau_i$ contains
$\tau_j$. Since $v(E)=\tau_j\cup\{\infty\}$ and
$\sharp(\tau_j)=q$, the determinant $D'$ of the matrix $\left(
\Upsilon_{I_s,k}\right)$, $1\le s\le q$, $k\in\tau_j$, does not
vanish at $\mu=m$. In order to prove the theorem it is enough to
show that there are $I_{q+1},\ldots,I_{q+\sharp_i}\in\mathcal N$
such that $v(I_s)=i$, $q+1\le s\le q+\sharp_i$, and the
determinant $D$ of the matrix $\left( \Upsilon_{I_s,k}\right)$,
$1\le s\le q+\sharp_i$, $k\in\tau_i$, does not vanish at $\mu=m$.
Set $I_{q+s+1}=(M-s,s,N-s)$, $M\le s\le N$, where $i=v(x^Mp^N)$.
By (\ref{e1}), (\ref{e2}) and (\ref{e3}),
\begin{equation}\label{ineq}
  g(\Upsilon_{I_s,k})<g(\Upsilon_{I_r,k}) \qquad \hbox{\rm  if }
  ~ k\ge i ~  \hbox{ \rm  and } ~ s\le q<r.
\end{equation}
Set $\lambda'=\left( \Upsilon_{I_s,k}\right)$, $q+1\le s\le
q+\sharp_i$, $k\in\tau_i\setminus\tau_j$. By lemma \ref{DDD},
$\det(\lambda'|_{\mu=m})\not=0$. Set
$\Upsilon_\varepsilon=\prod^{q+\sharp_i}_{s=1}\Upsilon_{I_s,\varepsilon
\gr{i}}$ for each bijection
$\varepsilon:\{1,\ldots,q+\sharp_i\}\to \tau_i$. By (\ref{ineq}),
$g(\Upsilon_\varepsilon)<g(D'\lambda'|_{\mu=m})$ if $\varepsilon
\gr{\{q+1,\ldots,q+\sharp_i\}}\not=\tau_i\setminus\tau_j$.
 Since
\[D'\lambda'|_{\mu=m}=\sum_{\varepsilon
\gr{\{q+1,\ldots,q+\sharp_i\}}=\tau_i\setminus\tau_j}
sign(\varepsilon)\Upsilon_\varepsilon,\]
  the product of the
leading monomials of $D'|_{\mu=m}$ and $\lambda'|_{\mu=m}$ is the
leading monomial of $D|_{\mu=m}$.
q.e.d.

\section{The moduli}

Set $s=s(n,m)=\inf (\Gamma(n,m)\backslash\langle n,m-n\rangle)$.
We say that (\ref{planecurve}) is in \em Legendrian short form \em
if $a_m=1$ and if $a_i\not=0$ and $i\in\Gamma(n,m)$, $i\in
\{m,s(n,m)\}$.

If $n=2$ or if $n=3$ and $m\in\{7,8\}$, $\Gamma(n,m)=\langle
n,m-n\rangle\supset\{m,\ldots\}$ and $x=t^n$, $y=t^m$ is the only
curve in Legendrian normal form such that the semigroup of its
conormal equals $\Gamma (n,m)$. If $n=3$ and $m\ge 10$ or if $n\ge
4$, $\langle n,n-m\rangle\not\supset\{m,\ldots,m+n-1\}$ and
$s(m,n)\in\{m,\ldots,m+n-1\}$.

\begin{lemma}
If \em (\ref{planecurve}) \em is in Legendrian normal form,
$\Gamma(n,m)\not=\langle n,m-n\rangle$ and the semigroup of the
conormal of \em (\ref{planecurve}) \em  equals $\Gamma(n,m)$,
$a_{s(n,m)} \neq 0$.
\end{lemma}
Proof.  Each $f\in\mathbb C\{x,y,p\}$ is congruent to a linear
combination of the series
\begin{equation}\label{ser}
    y,~nxp-my, ~x^i,~p^j, \qquad v(x^i),v(p^j)\le s
\end{equation}
modulo $(t^s)$. Since the series  (\ref{ser}) have different
valuations, one of these series must have valuation $s$,
$s\in\Gamma (n,m)\setminus \langle n,m-n\rangle$ and
$nxp-my=sa_st^s+\cdots$, $a_s\not=0$.
q.e.d.
Let $\mathcal{X}_{n,m}$ denote the set of plane curves
(\ref{planecurve}) such that (\ref{planecurve}) is in Legendrian
normal form and the semigroup of the conormal of
(\ref{planecurve}) equals $\Gamma(n,m)$. Let $W_n$ be the group of
$n$-roots of unity. There is an action of $W_n$ on
$\mathcal{X}_{n,m}$ that takes (\ref{planecurve}) into $x=t^n$,
$y=\sum_{i\ge m}\theta^{i-m}a_it^i$, for each $\theta\in W_n$. The
quotient $\mathcal{X}_{n,m}/W_n$ is an orbifold  of dimension
equal to the cardinality of the set $\{ m,...\} \backslash (
\Gamma(n,m)\setminus \{s(n,m)\})$.

\begin{theorem}
The set of isomorphism classes of generic Legendrian curves with
equisingularity type $(n,m)$ is isomorphic to
$\mathcal{X}_{n,m}/W_n$.
\end{theorem}
Proof. Let $\Lambda$ be a germ of an irreducible Legendrian
curve. There is a Legendrian map $\pi$ such that $\pi(\Lambda)$
has maximal contact with the curve $\{y=0\}$ and the tangent cone
of the conormal of $\Lambda$ equals $\{y=p=0\}$. Moreover, we can
assume that $\pi(\Lambda)$ has a parametrization of type
(\ref{planecurve}), with $a_m=1$. Assume that there is $i\in\Gamma
(m,n)$ such that $i\not=m,s(m,n)$ and $a_i\not=0$. Let $k$ be the
smallest integer $i$ verifying the previous condition. By lemmata
\ref{forget} and \ref{action} there are $a\in\mathbb C\{x,y,p\}$
and $\Phi \in \mathcal{J}$ such that $\imath^*a=a_kt^k+\cdots$ and
$\Phi$ takes   (\ref{planecurve}) into the plane curve $ x=s^n$,
$y=y(s)-a(s)+\delta$,  where $v(\delta) \geq 2 v (a)+m-2n$. Hence
we can assume that $a_i=0$ if $i\in\Gamma
(m,n)$, $i\not=m,s(m,n)$, and $i$ is smaller then 
the conductor $\sigma$ of the plane curve (\ref{planecurve}). 
There is a germ of diffeomorfism $\phi$ of the plane  that 
 takes the curve 
(\ref{planecurve}) into the curve 
$x=t^n$, $y=\sum_{i=m}^{\sigma-1} a_it^i$ (cf. \cite{zar}). This curve is 
in Legendrian normal form. The diffeomorphism $\phi$ induces an element of $\mathcal G$.

Let $\Phi$ be a contact transformation such
$\Phi(\mathcal{X})=\mathcal X$. Since the tangent cone of the
conormal of an element of $\mathcal{X}$ equals $\{y=p=0\}$, $\Phi
\in \mathcal{G}$. By theorem $\ref{subgroup}$,
$\Phi=\Psi\Psi_{\lambda,\mu}$, where $\Psi\in\mathcal{J}$ and
$\lambda,\mu\in\mathbb C^*$. Moreover, $\lambda\in W_n$ and
$\mu=\lambda^m$. By lemmata \ref{forget} and \ref{action},
$\Psi=Id$.
q.e.d.

\noindent
   A. Ara\'ujo and O. Neto \\
   CMAF, Universidade de Lisboa,\\
   Av. Gama Pinto, 2 1649, Lisboa\\
   Portugal\\
   {ant.arj@gmail.com}
   {orlando60@gmail.com}

\end{document}